\documentclass{amsart}
\usepackage{amssymb}
\newtheorem{theorem}{Theorem}

\newcommand{\QED}{\end{proof}}

\def\proclaim[#1]{{\bf #1}}
\def\BF#1.{{\bf #1.}}
\newcommand{\url}[1]{{\tt #1}}

\newcommand{\Godel}{G\"odel}

\newcommand{\N}{{\mathbb N}}

\newcommand{\R}{{\mathbb R}}

\newcommand{\set}[1]{\{\,{#1}\,\}}

\newcommand{\plus}{{+}}

\newcommand{\satisfies}{\models}

\newcommand{\Union}{\bigcup}

\newcommand{\smalllt}{\mathrel{\mathchoice{\raise2pt\hbox{$\scriptstyle<$}}{\raise1pt\hbox{$\scriptstyle<$}}{\raise0pt\hbox{$\scriptscriptstyle<$}}{\scriptscriptstyle<}}}
\newcommand{\smallleq}{\mathrel{\mathchoice{\raise2pt\hbox{$\scriptstyle\leq$}}{\raise1pt\hbox{$\scriptstyle\leq$}}{\raise1pt\hbox{$\scriptscriptstyle\leq$}}{\scriptscriptstyle\leq}}}
\newcommand{\lt}{\smalllt}

\newcommand{\boolval}[1]{\mathopen{\lbrack\!\lbrack}\,#1\,\mathclose{\rbrack\!\rbrack}}
\def\[#1]{\boolval{#1}}

\newcommand{\UnderTilde}[1]{{\setbox1=\hbox{$#1$}\baselineskip=0pt\vtop{\hbox{$#1$}\hbox to\wd1{\hfil$\sim$\hfil}}}{}}
\newcommand{\Undertilde}[1]{{\setbox1=\hbox{$#1$}\baselineskip=0pt\vtop{\hbox{$#1$}\hbox to\wd1{\hfil$\scriptstyle\sim$\hfil}}}{}}
\newcommand{\undertilde}[1]{{\setbox1=\hbox{$#1$}\baselineskip=0pt\vtop{\hbox{$#1$}\hbox to\wd1{\hfil$\scriptscriptstyle\sim$\hfil}}}{}}
\newcommand{\UnderdTilde}[1]{{\setbox1=\hbox{$#1$}\baselineskip=0pt\vtop{\hbox{$#1$}\hbox to\wd1{\hfil$\approx$\hfil}}}{}}
\newcommand{\Underdtilde}[1]{{\setbox1=\hbox{$#1$}\baselineskip=0pt\vtop{\hbox{$#1$}\hbox to\wd1{\hfil\scriptsize$\approx$\hfil}}}{}}

\newcommand{\st}{\mid}

\def\<#1>{\langle#1\rangle}

\newcommand{\ZFC}{{\rm ZFC}}
\newcommand{\ZF}{{\rm ZF}}

\newcommand{\CH}{{\rm CH}}

\newcommand{\AC}{{\rm AC}}

\newcommand{\AFA}{{\rm AFA}}
\newcommand{\AS}{{\rm AS}}

\newcommand{\MA}{{\rm MA}}

\newcommand{\MM}{{\rm MM}}

\newcommand{\PFA}{{\rm PFA}}

\newcommand{\HOD}{{\rm HOD}}

\newcommand{\cell}[1]{\boxit{\hbox to 17pt{\strut\hfil$#1$\hfil}}}
\newcommand{\head}[2]{\lower2pt\vbox{\hbox{\strut\footnotesize\it\hskip3pt#2}\boxit{\cell#1}}}
\newcommand{\boxit}[1]{\setbox4=\hbox{\kern2pt#1\kern2pt}\hbox{\vrule\vbox{\hrule\kern2pt\box4\kern2pt\hrule}\vrule}}
\newcommand{\Col}[3]{\hbox{\vbox{\baselineskip=0pt\parskip=0pt\cell#1\cell#2\cell#3}}}
\newcommand{\tapenames}{\raise 5pt\vbox to .7in{\hbox to .8in{\it\hfill input: \strut}\vfill\hbox to
.8in{\it\hfill scratch: \strut}\vfill\hbox to .8in{\it\hfill output: \strut}}}
\newcommand{\Head}[4]{\lower2pt\vbox{\hbox to25pt{\strut\footnotesize\it\hfill#4\hfill}\boxit{\Col#1#2#3}}}
\newcommand{\Dots}{\raise 5pt\vbox to .7in{\hbox{\ $\cdots$\strut}\vfill\hbox{\ $\cdots$\strut}\vfill\hbox{\
$\cdots$\strut}}}
\newcommand{\df}{\it} 
\begin{document}
\author{Joel David Hamkins}
\address{Department of Philosophy, New York University, 5
Washington Place, New York, NY 10003; and Mathematics, The
Graduate Center of The City University of New York, 365
Fifth Avenue, New York, NY 10016; and Mathematics, The
College of Staten Island of CUNY, Staten Island, NY 10314.}
\email{jhamkins@gc.cuny.edu, http://jdh.hamkins.org}
\today
\thanks{This article is based upon an argument I gave during the course of a three-lecture tutorial on set-theoretic geology at the summer school {\it Set Theory and Higher-Order Logic: Foundational Issues and Mathematical Developments}, at the University of London, Birkbeck in August 2011. Much of the article is adapted from and expands upon the corresponding section of material in \cite{Hamkins:TheSet-TheoreticalMultiverse}. My research has been supported in part by grants from the CUNY Research Foundation, the National Science Foundation (DMS-0800762) and the Simons Foundation, for which I am very grateful. I am grateful to Barbara Gail Montero for helpful comments.}
\begin{abstract}
The dream solution of the continuum hypothesis (\CH) would be a solution by which we settle the continuum hypothesis on the basis of a newly discovered fundamental principle of set theory, a missing axiom, widely regarded as true. Such a dream solution would indeed be a solution, since we would all accept the new axiom along with its consequences. In this article, however, I argue that such a dream solution to \CH\ is unattainable.
\end{abstract}

\title{Is the dream solution of the continuum hypothesis attainable?}\maketitle

Many set theorists yearn for a definitive solution of the continuum problem, what I call a {\df dream solution}, one by which we settle the continuum hypothesis (\CH) on the basis of a new fundamental principle of set theory, a missing axiom, widely regarded as true, which determines the truth value of \CH. In \cite{Hamkins:TheSet-TheoreticalMultiverse}, I describe the dream solution template as proceeding in two steps: first, one introduces the new set-theoretic principle, considered obviously true for sets in the same way that many mathematicians find the axiom of choice or the axiom of replacement to be true; and second, one proves the \CH\ or its negation from this new axiom and the other axioms of set theory. Such a situation would resemble Zermelo's proof of the ponderous well-order principle on the basis of the comparatively natural axiom of choice and the other Zermelo axioms. If achieved, a dream solution to the continuum problem would be remarkable, a cause for celebration.

In this article, however, I shall argue that a dream solution of \CH\ has become impossible to achieve. Specifically, what I claim is that our extensive experience in the set-theoretic worlds in which \CH\ is true and others in which \CH\ is false prevents us from looking upon any statement settling \CH\ as being obviously true. We simply have had too much experience by now with the contrary situation. Even if set theorists initially find a proposed new principle to be a natural, obvious truth, nevertheless once it is learned that the principle settles \CH, then this preliminary judgement will evaporate in the face of deep experience with the contrary, and set-theorists will look upon the statement merely as an intriguing generalization or curious formulation of \CH\ or $\neg\CH$, rather than as a new fundamental truth. In short, success in the second step of the dream solution will inevitably undermine success in the first step.

The continuum hypothesis, the assertion that every set of reals is either countable or equinumerous with the whole of $\R$, is equivalently formulated in \ZFC\ as the claim that the powerset $P(\N)$ of the natural numbers has the same cardinality as $\omega_1$, the first uncountable ordinal. Thus, the continuum hypothesis is the assertion that the two classical means of constructing uncountable sets give rise to the same uncountable cardinality. The question was open since the time of Cantor, appearing at the top of Hilbert's famous 1900 list of open questions, until \Godel\ proved that $\ZFC+\CH$ holds in the constructible universe $L$ of any model of \ZFC, and Cohen proved that $L$ has a forcing extension $L[G]$ satisfying $\ZFC+\neg\CH$. Going beyond this connection to the constructible universe, both the continuum hypothesis and its negation are forceable over any model of set theory:

\begin{theorem}[Cohen, Solovay]\label{Theorem.CHnotCH} The set-theoretic universe $V$ has forcing extensions
 \begin{enumerate}
 \item $V[G]$, collapsing no cardinals, such that $V[G]\satisfies\neg\CH$, and
 \item $V[H]$, adding no new reals, such that $V[H]\satisfies\CH$.
 \end{enumerate}
\end{theorem}

Although it was formerly common to undertake forcing constructions only for countable transitive models of fragments of \ZFC, one may formalize the forcing method as an internal \ZFC\ construction, rather than a meta-theoretic construction, and thereby make sense of forcing over an arbitrary model of set theory, including the over the full universe $V$. By means of what I have called the naturalist account of forcing \cite{Hamkins:TheSet-TheoreticalMultiverse}, or by the classical Boolean ultrapower (Boolean-valued quotient) approach (e.g. see \cite{HamkinsSeabold:BooleanUltrapowers, FuchsHamkinsReitz:Set-theoreticGeology}), one may legitimize the forcing-over-the-universe approach to forcing, by now the most common approach in the set-theoretic literature, and it is this approach that resonates most strongly with the multiverse perspective that the forcing extensions of $V$ are real.

The theorem shows that every model of set theory is very close to models with the opposite answer to \CH. Since the $\CH$ and $\neg\CH$ are easily forceable, the continuum hypothesis is something like a lightswitch, which can be turned on and off by moving to ever larger forcing extensions. Indeed, a key concept in the modal logic of forcing
\cite{HamkinsLoewe2008:TheModalLogicOfForcing} is that of a {\df switch}, a statement $\psi$ of set theory such that both $\psi$ and $\neg\psi$ are forceable over any forcing extension of the universe. Thus, the continuum hypothesis is a switch, one which set theorists today have truly flicked many times. In each
case of the theorem the forcing is relatively mild, with the new universes, for example, having all the same large cardinals as the original universe, a fact that refuted \Godel's hope that large cardinals might settle the \CH. After decades of experience and study, set-theorists today have a profound understanding of how to achieve the continuum hypothesis or its negation in diverse models of set theory---forcing it or its negation in innumerable ways, while simultaneously controlling other set-theoretic properties---and have therefore come to a deep knowledge of the extent of the continuum hypothesis and its negation in set-theoretic worlds.

More generally, I have argued in \cite{Hamkins:TheSet-TheoreticalMultiverse} that in set theory we have come to discover an entire multiverse of set concepts, each giving rise to a corresponding set-theoretic universe. Many of these set-theoretic universes are the universes that we have long known and described in set theory, such as $L$, $\HOD$, $L[\mu]$, $K$, $\HOD[A]$, $V_\kappa$, $H_{\kappa^\plus}$, $L(\R)$, $L(V_{\lambda+1})$ and innumerably many others, including especially the forcing extensions of these and other models by any of the enormous collection of forcing notions that have been studied in set theory. In the past half-century, set theorists have gained a precise familiarity with the nature of these diverse set-theoretic worlds; we move from one to another with ease. Part of my goal in the multiverse article was to tease apart two often-blurred aspects of set-theoretic Platonism, namely, to separate the claim that the set-theoretic universe has a real mathematical existence from the claim that it is unique. The multiverse perspective is meant to affirm the realist position, while denying the uniqueness of our set-theoretic background concept. What we have learned in set theory is that we have a choice of diverse set-theoretic universes, each arising from its own iterated set concept. We may regard each of these universes to be fully as real as the universists take their background set concept to be, just as one may regard all the shades of blue as actual colors, regardless of any debate about which of them is to be deemed officially ``blue''. Our mathematical experience is that these alternative set-theoretic worlds are perfectly fine set-theoretically; there is absolutely nothing wrong with them, and we have gained an extensive experience living in them.

On the multiverse view the continuum hypothesis is a settled question, for the answer consists of the expansive, detailed knowledge set theorists have gained about the extent to which the \CH\ holds and fails in the multiverse, about how to achieve it or its negation in combination with other diverse set-theoretic properties. Of course, there are and will always remain questions about whether one can achieve \CH\ or its negation with this or that hypothesis, but the point is that the most important and essential facts about \CH\ are deeply understood, and it is these facts that constitute the answer to the \CH\ question.

One way to understand my argument that we will not achieve the dream solution of the continuum problem is to phrase it in terms of Quine's web of belief \cite{Quine1951:TwoDogmasOfEmpricism}, the view that our system of beliefs, even our scientific beliefs, forms a coherent interconnected web for which in principle any part could be altered in the light of new knowledge, but also for which new evidence needn't force the revision of any particular given part. My argument is that our extensive experience in all the various set-theoretic universes with \CH\ and with $\neg\CH$, especially those obtained by forcing, have pushed the dual possibilities of \CH\ and $\neg\CH$ to the center of our set-theoretic beliefs. Specifically, the idea that in any set-theoretic context, there is a nearby and perfectly legitimate set-theoretic universe where \CH\ has a prescribed value has become a core set-theoretic intuition at the center of our beliefs about \CH\ in set theory, and set theorists would tenaciously hold onto it come what may. Any newly proposed set-theoretic principle settling \CH, if regarded as a new fundamental truth, therefore, could not consistently be incorporated into the web of belief without overturning these central commitments, but if regarded merely as an intriguing formulation or strengthening of \CH\ or $\neg\CH$, in contrast, could be incorporated into the web of belief with ease.

To support this claim further, let me borrow from the opposition. In his excellent recent article, Daniel Isaacson \cite{Isaacson2008:TheRealityOfMathematicsAndTheCaseOfSetTheory} mounts a vigorous appeal to second-order categoricity arguments in order to establish the uniqueness of the set-theoretic universe and thereby establish the cumulative hierarchy as what he calls a particular mathematical structure. This style of argument goes back to Zermelo's categoricity argument characterizing the models of second-order set theory as precisely those of the form $V_\kappa$ for an inaccessible cardinal $\kappa$, but Isaacson proceeds further to provide a characterization of the entire set-theoretic universe in terms of second-order set-theoretic properties. Isaacson's argument shares important features with Martin's similar argument of \cite{Martin2001:MultipleUniversesOfSetsAndIndeterminateTruthValues}, and both Isaacson and Martin eventually use the categoricity arguments in order to defend the view that the continuum hypothesis is a definite mathematical question (though neither says which way the answer goes). Thus, Isaacson is aligned with the universist position. Nevertheless, let us look more closely at his description of how particular mathematical structures become established, of how we come to know them, to support my argument against the universe view.

Specifically, Isaacson distinguishes between the {\df particular} as opposed to {\df general} mathematical structures, discussing the distinction at length in its connection to structuralism. The structures of the natural numbers, the integers and the real numbers are each particular mathematical structures, as opposed to the class of all groups, all rings or all topological spaces, which are examples of classes of general mathematical structures. Although both particular and general mathematical structures are treated by means of mathematical axioms, the nature of this treatment is fundamentally different in the two cases. Namely, with the particular structures, we identify general principles true of these structures, which encapsulate our knowledge about them, characterize them up to isomorphism, and serve as axioms in the sense that we use those general principles to unify diverse arguments about the structures. With general mathematical structures, in contrast, such as when we specify the class of rings by the ring axioms, the axioms have the character of definitions or demarcations of the domain of discourse rather than self-evident truths or encapsulations of specific knowledge.

Isaacson explains that our knowledge of particular mathematical structures arises from our informal mathematical practice and experience with them.

\begin{quote}
``If the mathematical community at some stage in the development of mathematics has succeeded in becoming (informally) clear about a particular mathematical structure, this clarity can be made mathematically exact. Of course by the general theorems that establish first-order languages as incapable of characterizing infinite structures the mathematical specification of the structure about which we are clear will be in a higher-order language, usually by means of a full second-order language. Why must there {\it be} such a characterisation? Answer: if the clarity is genuine, there must be a way to articulate it precisely. If there is no such way, the seeming clarity must be illusory. Such a claim is of the character as the Church-Turing thesis, for every apparently algorithmic process, there is a Turing machine or $\lambda$-calculus formal computation. In the present case, for every particular structure developed in the practice of mathematics, there is [a] categorical characterization of it.''(p. 31, December 20, 2007 version)
\end{quote}

\noindent He specifically identifies a notion of informal rigour by which we come to understand the particular mathematical structures.

\begin{quote}
``The something more that is needed to represent true sentence of arithmetic as logical consequences of the second-order axioms of arithmetic is\ldots the informal rigour by which we have come to understand these second-order axioms, and thereby to see that they are coherent. It is a development of mathematical understanding through informal rigour and not some further derivation that is needed.\ldots We must reflect on our conceptual understanding of a given particular mathematical structure as it has developed to see how it is that truths of e.g. arithmetic are those that hold in the structure of the natural numbers which we have succeeded in characterizing.'' (p. 33, December 20, 2007 version)
\end{quote}

I claim that it is precisely this kind of mathematical experience that set-theorists have gained with respect to the various (particular) forcing extensions of the universe and the other models of set theory commonly studied. We know what it is like to live in the universe obtained by adding $\aleph_2$ many Cohen reals over $L$ and what it is like in the $L(\R)$ of that extension. We know what it is like in the universe obtained by forcing $\MA+\neg\CH$ over $L$, or by forcing with Sacks forcing, or by iteratively adding a dominating real. The resulting universes are places we've been. We are deeply familiar with the universe obtained by forcing with the Laver preparation of a supercompact cardinal and with the model of \PFA\ obtained by Baumgartner's similar construction. This is not to say that we know everything that there is to know about these universes---there will always be more to learn about them---just as Isaacson mentions that we do not know the complete theory of arithmetic even though we seem to have the integers as a particular mathematical structure. But we have sufficient experience living in these diverse set-theoretic worlds to know that the concept of set used in each of them is perfectly robust and satisfactory as a concept of set. Each of these universes feels fully set-theoretic, and one can imagine living out a full mathematical life inside almost any one of them.

It is for this reason that the dream solution has become impossible. Our situation with \CH\ is not merely that \CH\ is formally independent and we have no additional knowledge about whether it is true or not. Rather, we have an informed, deep understanding of how it could be that \CH\ is true and how it could be that \CH\ fails. We know how to build the \CH\ and $\neg\CH$ worlds from one another. Set theorists today grew up in these worlds, comparing them and moving from one to another while controlling other subtle features about them. Consequently, if someone were to present a new set-theoretic principle $\Phi$ and prove that it implies $\neg\CH$, say, then we could no longer look upon $\Phi$ as obviously true. To do so would negate our experience in the \CH\ worlds, which we found to be perfectly set-theoretic. It would be like someone proposing a principle implying that only Brooklyn really exists, whereas we already know about Manhattan and the other boroughs. And similarly if $\Phi$ were to imply $\CH$. We are simply too familiar with universes exhibiting both sides of \CH\ for us ever to accept as obviously true a principle that is false in some of them. So success in the second step of the dream solution fatally undermines success in the first step.

In summary, for any particular attempt at a dream solution of the continuum problem, where a new set-theoretic principle is proposed and proved to settle \CH, then I predict that set theorists will object to the claim that the principle is obviously true, and furthermore, their objections will arise from a wellspring of deep mathematical experience with the contrary hypothesis.

Let me illustrate with examples of how this predicted pattern of response has in fact already occurred a few times. Consider first the reaction to Chris Freiling's delightful axiom of symmetry \cite{Freiling1986:AxiomsOfSymmetry:ThrowingDartsAtTheRealLine}. Freiling describes his argument as containing ``a simple philosophical `proof' of the negation of Cantor's continuum hypothesis'' and presents a line of reasoning that I find exactly to follow the dream solution template. Namely, Freiling begins with a bit of philosophical intuition-building, ``subject[ing] the continuum to certain thought experiments involving random darts,'' before ultimately landing at his axiom of symmetry, presented as a true and natural axiom, an ``intuitively clear axiom,'' whose truth, Freiling argues, follows from our strongest pre-reflective ideas about symmetry and likelihood, the same intuitions that underlie the fundamental concepts of measure theory.

Let us jump into the details, which are fun to consider. The axiom of symmetry is the assertion that for any function $f$ mapping reals to countable sets of reals, there are real numbers $x$ and $y$ such that $y\notin f(x)$ and $x\notin f(y)$. To argue for the natural appeal of the axiom, Freiling proposes that we imagine throwing two darts in succession at a dart board, considering precisely where they land. The first dart lands at some position $x$, and because $f(x)$ is a countable set, we expect almost surely that the second dart will land at a point $y$ not in that set, and so almost surely $y\notin f(x)$. But since the order in which we consider the darts shouldn't seem to matter, we conclude by symmetry that almost surely $x\notin f(y)$ as well. So almost surely our darts will land at positions $x$ and $y$ fulfilling the axiom of symmetry claim $y\notin f(x)$ and $x\notin f(y)$. Freiling emphasizes that not only do we have natural reason to expect that there is a pair $(x,y)$ with the desired property, but what is more, we should expect that almost all pairs have the desired property: ``actually [the axiom], being weaker than our intuition, does not say that the two darts have to do anything. All it claims is that what heuristically will happen every time, can happen.'' Thus, Freiling argues for the natural appeal and truth of the axiom of symmetry.

Many mathematicians, in my personal experience, though usually not the more experienced mathematicians, find this dart-throwing thought experiment quite appealing upon the first presentation, at least at this point in the argument, and regard it as providing a natural support for the axiom of symmetry. These mathematicians are often quite surprised, in the second part of Freiling's argument, to learn that the axiom of symmetry is actually equivalent to $\neg\CH$, and they universally revise their initial judgement of the axiom as a result. Nevertheless, the equivalence of \AS\ with $\neg\CH$ is not difficult to establish. The forward implication is easy, for if \CH\ holds, then there is a well-ordering of $\R$ in order type $\omega_1$, and we may consider the function $f$ mapping every real $x$ to the initial segment of the order up to $x$, a countable set. The point is that for any two real numbers, either $x$ precedes $y$ or conversely in the well-order, and so either $x\in f(y)$ or $y\in f(x)$, contrary to the axiom of symmetry. So $\AS\implies\neg\CH$, and this is already enough for the intended dream solution. Conversely, if $\CH$ fails, then for any choice of $\omega_1$ many distinct reals $x_\alpha$, for $\alpha<\omega_1$, the union $\Union_{\alpha<\omega_1} f(x_\alpha)$ has size $\omega_1$ and so by $\neg\CH$ there must be a real $y$ not in any $f(x_\alpha)$. Since $f(y)$ contains at most countably many $x_\alpha$, there must be some $x_\alpha$ with $x_\alpha\notin f(y)$ and since we've already ensured $y\notin f(x_\alpha)$, we have the desired pair (and in fact many such pairs). In summary, the axiom of symmetry is equivalent to $\neg\CH$, and Freiling has exactly carried out the dream solution template for \CH.

Was Freiling's argument received as the longed-for solution to \CH? No. Many mathematicians objected that Freiling's argument  implicitly assumes for a given function $f$ that various sets are measurable, including most importantly the set $\set{(x,y)\st y\in f(x)}$ in the plane. Although each vertical slice of this set is countable and hence measure zero, Freiling's argument relies on our intuitions concerning two independent events---the two dart throws---and depends fundamentally on the robustness of our measure concepts for this two-dimensional set. Freiling anticipates this objection, answering that the intuitive justification he offers for the axiom of symmetry can be viewed as prior to the mathematical development of measure theory. Our confidence in the axiom of symmetry itself, he argues, rests on the same philosophical or pre-reflective ideas that underlie the technical mathematical requirements we impose on our measure theory in the first place, such as our insistence that measures be countably additive, and we would therefore seem to have as much direct support for the axiom of symmetry as we have for those requirements.

My main point here is to observe the nature of the criticism, rather than to debate the merits of the reply. Specifically, I want to call attention to the fact that mathematicians objected to Freiling's argument largely from a perspective of deep experience and familiarity with non-measurable sets and functions, including extreme violations of the Fubini property of the kind on which his argument relies. For mathematicians with this experience and familiarity, the pre-reflective arguments simply fall flat. We have become deeply skeptical of any intuitive or naive use of measure concepts precisely because we know the pitfalls. Because of our mathematical experience, we know how complicated and badly-behaved functions and sets of reals can be in their measure-theoretic properties. We know that any naive account of measure will have a fundamental problem dealing with subsets of the plane all of whose vertical sections are countable and all of whose horizontal sections are co-countable, for example, precisely because the sets looks very small from one direction and very large from another direction, while we expect that rotating a set should not change its size.

Consider the variation of Freiling's argument based directly on that intuition: intuitively, all sets are measurable, and also rotating sets in the plane preserves measure; but if the \CH\ were to hold, then there is a subset of the unit square, such as the graph of a well-ordering of the unit interval in order type $\omega_1$, with all horizontal sections countable and all vertical sections co-countable; this set appears to have measure zero from one direction and measure one from another. Hence, \CH\ fails.

The opponents to Freiling's argument will answer this modified argument by pointing out that our detailed experience with non-measurable sets prevents us from accepting the naive claim that every set is measurable. Similarly, for the original Freiling argument, we are simply not convinced by Freiling's argument that the axiom of symmetry is intuitively true, even if he is using the same intuitions that guided us to the basic principles of measure theory. An experienced mathematician answers Freiling's intuitive appeal by pointing out that it relies fundamentally on having nice measure-theoretic properties for the set $\set{(x,y)\st y\in f(x)}$, whereas we have extensive experience with very badly behaved sets and functions, and no reason to suppose this set is not also badly behaved in that way. In an extreme instance of this, inverting Freiling's argument, set theorists sometimes reject the axiom of symmetry as a fundamental axiom, precisely because of the counterexamples to it that one can produce under \CH. In short, because we are deeply familiar with a way that the axiom of symmetry can fail, we do not accept the intuitive justifications as establishing it as true.

Ultimately, rather than being accepted as the longed-for solution to the continuum hypothesis, Freiling's argument is instead most often described as providing an attractive equivalent formulation of $\neg\CH$, a curious and interesting form of it. The typical presentation of the axiom of symmetry includes a discussion of Freiling's dart-throwing justification, but in my experience this discussion is usually given not as evidence that the axiom is true, but rather as a warning about the measure-theoretic monsters, a warning that we must take extra care with issues of non-measurability, lest we be fooled. In this way, Freiling's simple philosophical argument is turned on its head, used not as a justification of the axiom, but rather as a warning about the error that may arise from a naive treatment of measure concepts, a warning that what seems obviously true might still be wrong. The entire episode bears out the pattern of response I predict for any attempted use of the dream solution template, namely, a rejection of the new axiom from a perspective of deep mathematical experience with the contrary.

Let me turn now to a second illustration of this pattern of response. Consider the set-theoretic principle that I have called the {\df powerset size axiom} {\rm PSA}, the axiom asserting plainly that smaller sets have fewer subsets:
$$\forall x,y\qquad |x|\lt |y|\implies |P(x)|\lt |P(y)|.$$
Set theorists understand the situation of this axiom very well, which I shall shortly explain. But how is it received in mathematics generally? Extremely well! A large number of mathematicians, including some very good ones, although invariably from non-logic-related areas of mathematics, look favorably upon the axiom when it is first considered, viewing it as highly natural or even obviously true. They take the axiom to express what seems be a basic intuitive principle, namely, that a strictly smaller set should have strictly fewer subsets. The principle, for example, is currently the top-rated answer \cite{HamkinsMO6594:AnswerToWhatAreSomeReasonable-soundingStatementsIndepOfZFC?} among dozens to a popular mathoverflow question seeking examples of reasonable-sounding statements that are nevertheless independent of the axioms of set theory, and the same issue has arisen in at least three other mathoverflow questions, posted by mathematicians asking naively whether the PSA is true, or how to prove it or indeed asking with credulity how it could not be provable. My experience is that a brief conversation with mathematicians at your favorite math tea stands a good chance to turn up additional examples of mathematicians who find the axiom to express a basic fact about sets.

Meanwhile, set theorists almost never agree with this assessment. They know that the axiom is independent of \ZFC, for one can achieve all kinds of
crazy patterns for the continuum function $\kappa\mapsto 2^\kappa$ via Easton's theorem. Even Cohen's original model of $\neg\CH$ had $2^\omega=2^{\omega_1}$, the assertion known as Luzin's hypothesis \cite{Luzin1935:SurLesEnssemblesAnalytiquesNuls}, which had been proposed as an alternative to the continuum hypothesis. Furthermore, Martin's axiom implies $2^\omega=2^\kappa$ for all $\kappa<2^\omega$, which can mean additional violations of $\text{PSA}$ when \CH\ fails badly. So not only do set-theorists know that $\text{PSA}$ can fail, but also they know that $\text{PSA}$ must fail in models of the axioms, such as the proper forcing axiom \PFA\ or Martin's maximum \MM, that are often favored particularly by set-theorists. For some set theorists, the simple philosophical support of $\text{PSA}$ is suppressed as naive in favor of a complex philosophical interest in the forcing axioms, which imply $\neg\text{PSA}$.

So the situation with the powerset size axiom is that a set-theoretic principle (1) that many mathematicians find to be obviously true, (2) which surely expresses an intuitively clear pre-reflective principle about the concept of size, and (3) which furthermore is known by set-theorists to be perfectly safe in the sense that it is relatively consistent with the other axioms of \ZFC\ and in fact a consequence of the generalized continuum hypothesis, is nevertheless almost universally rejected by set-theorists when it is proposed as a fundamental axiom.

Although the $\text{PSA}$ does not settle the \CH, my point here is that nevertheless this rejection of the powerset size axiom follows the predicted pattern of response to the dream solution template: a clear and succinct mathematical principle, proposed as fundamental, enjoys a strong intuitive appeal and obvious nature, but is nevertheless rejected from a perspective of deep experience with the contrary. We simply know too much about the various ways that the principle can be violated and have too much experience working in models of set theory where the principle fails to accept it as a fundamental truth.

Let me propose a thought experiment of my own. Imagine that the history of set theory had proceeded differently, that the powerset size axiom had been considered at the very beginning of set theory, perhaps used in a proof settling a major open question of the period. For example, consider the question whether the symmetric groups $S_\kappa$ and $S_\lambda$ on distinct infinite cardinals, the respective groups of all permutations on fixed sets of size $\kappa$ and $\lambda$, are non-isomorphic as groups. The PSA implies an immediate affirmative answer, since the groups have size $2^\kappa$ and $2^\lambda$, respectively, and under the PSA these are distinct.\footnote{Meanwhile, in fact, no additional axiom is needed to settle the question, even if it should happen that $2^\kappa=2^\lambda$, for the groups $S_\kappa$ and $S_\lambda$ can be distinguished by the Schreier--Ulam--Baer theorem.} Because of the natural appeal of the axiom, it seems plausible to imagine that the PSA might have found its way onto the standard list of axioms: perhaps an alternative Zermelo might have formalized the symmetric group argument and presented $\ZFC+\text{PSA}$ as the list of fundamental axioms, just as our actual Zermelo formalized his proof of the well-order principle with ZC. In this imaginary alternative history, we might now look upon models of $\neg\text{PSA}$ as strange in some fundamental way, regarding them to violate a basic intuitive principle of sets concerning the relative sizes of power sets; perhaps our reaction to these models would be like the current reaction some mathematicians (not all) have to models of $\ZF+\neg\AC$ or to models of Aczel's anti-foundation axiom \AFA, namely, the view that the models may be interesting mathematically and useful for a purpose, but ultimately they violate a basic principle of sets. The point I want to make with this thought experiment is that without our current detailed technical knowledge of how PSA can fail, we would likely have found the intuitive appeal more compelling. Indeed, the ease with which set theorists today shrug off the enormous intuitive pull of PSA is surprising.

I have argued, then, that there will be no dream solution of the continuum hypothesis. Let me now go somewhat beyond this claim and issue a challenge to those who propose to solve the continuum problem by some other means. My challenge to anyone who proposes to give a particular, definite answer to \CH\ is that they must not only argue for their preferred answer, mustering whatever philosophical or intuitive support for their answer as they can, but also they must explain away the illusion of our experience with the contrary hypothesis. Only by doing so will they overcome the response I have described, rejection of the argument from extensive experience of the contrary. Before we will be able to accept \CH\ as true, we must come to know that our experience of the $\neg\CH$ worlds was somehow flawed; we must come to see our experience in those lands as illusory. It is insufficient to present a beautiful landscape, a shining city on a hill, for we are widely traveled and know that it is not the only one.

\bibliographystyle{alpha}
\bibliography{MathBiblio,HamkinsBiblio}

\end{document}